\documentclass[11pt,reqno]{amsart}

\usepackage{graphics}
\usepackage{amssymb, a4wide}
\usepackage{amsmath}
\usepackage{graphicx}
\usepackage{epsfig}

\newtheorem{proposition}{Proposition}
\newtheorem{theorem}{Theorem}
\newtheorem{example}{Example}
\newtheorem{remark}{Remark}
\newtheorem{lemma}{Lemma}
\newtheorem{corollary}{Corollary}

\title{On space-like constant slope surfaces and Bertrand curves in Minkowski 3-space}

\author{Murat Babaarslan}
\address[M. Babaarslan]{Bozok University,
Department of Mathematics,
66100, Yozgat, Turkey}
\email{murat.babaarslan@bozok.edu.tr}

\author{Yusuf Yayli}
\address[Y. Yayli]{Ankara University,
Department of Mathematics,
06100, Ankara, Turkey}
\email{yayli@science.ankara.edu.tr}

\subjclass[2010]{53B25}

\keywords{Bertrand curve, helix, pseudo-spherical Darboux image, Lorentzian Sabban frame, pseudo-spherical evolutes, space-like constant slope surface.}

\date{\today}

\begin{document}

\parindent 0mm
\parskip 2mm

\maketitle

\begin{abstract}

In the present paper, we define the notions of Lorentzian Sabban frames and de Sitter evolutes of the unit speed space-like curves on de Sitter 2-space $\mathbb{S}^{2}_{1}$. In addition, we investigate the invariants and geometric properties of these curves. Afterwards, we show that space-like Bertrand curves and time-like Bertrand curves can be constructed from unit speed space-like curves on de Sitter 2-space $\mathbb{S}^{2}_{1}$ and hyperbolic space $\mathbb{H}^{2}$, respectively. We obtain the relations between Bertrand curves and helices. Also we show that pseudo-spherical Darboux images of Bertrand curves are equal to pseudo-spherical evolutes in Minkowski 3-space $\mathbb{R}^{3}_{1}$. Moreover we investigate the relations between Bertrand curves and space-like constant slope surfaces in $\mathbb{R}^{3}_{1}$. Finally, we give some examples to illustrate our main results.

\end{abstract}


\section{Introduction}

A lot of well-known geometric objects are defined with respect to making a constant angle with a given, distinguished direction. Firstly, helices are curves which make a constant angle with a fixed direction \cite{garnica}. There are a lot of interesting applications of helices, for example; $\alpha$-helices, DNA double and collagen triple helix, nano-springs, carbon nano-tubes, K-helices, helical staircases, helical structures in fractal geometry and so on \cite{munteanu}. A second example is logarithmic spirals or spira mirabilis extensively investigated by Jacob Bernoulli which make a constant angle with the radial direction. They are very fascinating curves since they are self-similar.  Possibly as a result of this unique property, they appear in nature, including in shells, sunflower heads and horns (see \cite{boyadzhiev}). A third famous example is loxodromes also known as rhumb lines which cut all of the meridians at a constant angle on the Earth and thus they are important in navigation \cite{garnica}. For their interesting applications, we refer to \cite{babaarslan2, babaarslan3} and references therein.

Another interesting and important notion is evolute curve. Evolute curve was discovered by C. Huygens in 1673 while trying to develop the accurate pendulum clock. The curve $\alpha$ is called as evolute of $\tilde{\alpha}$ and the curve $\tilde{\alpha}$ is called as involute of $\alpha$ if the tangent vectors are perpendicular at the corresponding points for any $s \in I\subset \mathbb{R}$. Also, the locus of the center of the curvatures of $\alpha$ is evolute or focal curve of $\alpha$ \cite{millman}.

In addition, Bertrand curves discovered by J. Bertrand in 1850 are defined as a special curve which shares its principal normals with another special curve (called Bertrand mate). We say that the curve $\alpha$ is Bertrand curve if and only if there are non-zero real numbers $A, B$ such that $A\kappa(s)+B\tau(s)=1$  for any $s\in I\subset \mathbb{R}$. Thus, plane curves and circular helices are Bertrand curves. Also, if $\beta$ and $\gamma$ are two different involutes of $\alpha$, then they are Bertrand mates. Bertrand curves are particular examples of parallel (offset) curves which are used in computer-aided design (CAD) and computer-aided manufacture (CAM) \cite{nutbourne}.

The fact that Bertrand curves can be constructed from unit speed curves on 2-sphere $\mathbb{S}^{2}$ was shown by Izumiya and Takeuchi in \cite{izumiya1}. They defined spherical evolutes of the unit speed curves on $\mathbb{S}^{2}$ and showed that these spherical evolutes are locus of the center of the curvatures of unit speed curves on $\mathbb{S}^{2}$. Also, they proved that spherical Darboux images of Bertrand curves are equal to spherical evolutes of unit speed curves on $\mathbb{S}^{2}$ in Euclidean 3-space $\mathbb{R}^{3}$.

After that, Izumiya et al. \cite{izumiya2} introduced hyperbolic Frenet-Serret formulae of unit speed space-like curves on hyperbolic space $\mathbb{H}^{2}$ and defined hyperbolic evolutes of these curves. Also they studied geometric properties of them. They showed that these hyperbolic evolutes are the locus of the center of geodesic curvatures of the unit speed space-like curves on $\mathbb{H}^{2}$ in Minkowski 3-space $\mathbb{R}^{3}_{1}$.

As a generalization of helices and logaritmic spirals, we can think constant slope surfaces \cite{munteanu} making a constant angle with the position vectors. Munteanu \cite{munteanu} gave the characterization of these surfaces and obtained some beautiful shapes in $\mathbb{R}^{3}$.

In \cite{babaarslan1}, we found a parametrization of Bertrand curves different from in \cite{izumiya1} and gave some relations between Bertrand curves and constant slope surfaces in $\mathbb{R}^{3}$.

Another study in this direction is \cite{fu}, where Fu and Yang gave a complete classification of space-like constant slope surfaces in $\mathbb{R}^{3}_{1}$. As it was mentioned before, the position vectors of these surfaces make a constant angle with the normal at each point on the surfaces. Fu and Yang [5] showed that the surface $S \subset \mathbb{R}^{3}_{1}$ is a space-like constant slope surface lying in the space-like cone if and only if it can be parametrized by
\begin{equation}
\label{eq1}
x(u,v)=u\cosh\theta \big(\cosh\xi_{1}f(v)+\sinh\xi_{1}f(v)\times f'(v) \big),
\end{equation}
where $\theta$ is a constant different from $0$, $\xi_{1}=\xi_{1}(u)=\tanh \theta \ln u$ and $f$ is a unit speed space-like curve on de Sitter 2-space $\mathbb{S}^{2}_{1}$.

Also, the surface $S \subset \mathbb{R}^{3}_{1}$ is a space-like constant slope surface lying in the time-like cone if and only if it can be parametrized by
\begin{equation}
\label{eq2}
x(u,v)=u\sinh\theta \big(\cosh\xi_{2}g(v)+\sinh\xi_{2}g(v)\times g'(v) \big),
\end{equation}
where $\theta$ is a constant different from $0$, $\xi_{2}=\xi_{2}(u)=\coth \theta \ln u$ and $g$ is a unit speed space-like curve on hyperbolic space $\mathbb{H}^{2}$.

Constant slope surfaces have beautiful shapes and they are interesting in terms of differential geometry. The study of these surfaces is similar to that of helices and logarithmic spirals. At least for their shapes, we can say that constant slope surfaces are one of the most fascinated surfaces in both Euclidean and Minkowski 3-space \cite{munteanu}.

In this paper,  we will define the notions of Lorentzian Sabban frames and de Sitter evolutes of the unit speed space-like curves on de Sitter 2-space $\mathbb{S}^{2}_{1}$. Also we will investigate the invariants and geometric properties of these curves. Afterwards, we will show that space-like Bertrand curves and time-like Bertrand curves can be constructed from unit speed space-like curves on de Sitter 2-space $\mathbb{S}^{2}_{1}$ and hyperbolic space $\mathbb{H}^{2}$, respectively. We will obtain some relations between Bertrand curves and helices. Also we will show that pseudo-spherical Darboux images of Bertrand curves are equal to pseudo-spherical evolutes in $\mathbb{R}^{3}_{1}$. Moreover we will investigate some relations between Bertrand curves and space-like constant slope surfaces in $\mathbb{R}^{3}_{1}$. Finally, we will give some examples of our main results and the corresponding pictures will be drawn via Mathematica.


\section{Basic notations, definitions and formulas}

In this section, we give the basic notations and some results in Minkowski 3-space. For more details, we refer to \cite{lopez, o'neill}. We call the real vector space $\mathbb{R}^{3}$ endowed with the standard Lorentzian metric
\begin{eqnarray}
\label{eq1}
\langle x,y\rangle=x_1y_1+x_2y_2-x_3y_3 \nonumber
\end{eqnarray}
as Minkowski 3-space $\mathbb{R}^{3}_{1}$, where $x=(x_1,x_2,x_3)$ and  $y=(y_1,y_2,y_3)$ are vectors in $\mathbb{R}^3$. We say that a vector $x\in \mathbb{R}^{3}_{1}$ is space-like if $\langle x,x\rangle>0$ or  $x=0$, it is time-like if $\langle x,x\rangle<0$ and it is light-like (null) if $\langle x,x\rangle=0$ and  $x\neq0$. This category of the given vector $x\in \mathbb{R}^{3}_{1}$ is called its causal character. The pseudo norm (length) of a vector  $x$ is given by  $\|x\|=\sqrt{|\langle x,x\rangle|}$.

Given a regular curve $\alpha:I\subset\mathbb{R}\rightarrow\mathbb{R}^{3}_{1}$, we say that $\alpha$ is space-like (resp. time-like, light-like) if all of its velocity vectors $\alpha'(t)$ are space-like (resp. time-like, light-like).

If $\alpha$ is a space-like or time-like curve, then we say that $\alpha$ is a non-light-like curve. In this case, there exist a change of parameter $t$, that is, $s(t)$ such that $\|\alpha'(s)\|=1$. Then we say that $\alpha$ is parametrized by the arc-length parameter. In this case, $\alpha$ is called as a unit speed curve.

In Minkowski 3-space $\mathbb{R}^{3}_{1}$, the Lorentzian cross-product of any vectors $x=(x_1,x_2,x_3) \ \textrm{and}\ y=(y_1,y_2,y_3)\in \mathbb{R}^{3}_{1}$ is defined as follows:
\begin{eqnarray}
x\times y=\left|%
\begin{array}{ccc}
  e_{1} & e_{2} & -e_{3}\\
  x_1 & x_2 & x_3\\
  y_1 & y_2 & y_3\\
\end{array}%
\right|=(x_2y_3-x_3y_2,x_3y_1-x_1y_3,x_2y_1-x_1y_2).\nonumber
\end{eqnarray}
Like the cross-product in $\mathbb{R}^{3}$, the Lorentzian cross-product has similar algebraic and geometric properties, for example:
\begin{enumerate}
\item[ i.] $\langle x\times y,z\rangle= \det(x,y,z)$;\\
\item[ii.] $x\times y=-y\times x$;\\
\item[iii.] $(x\times y)\times z=-\langle x,z\rangle y+\langle y,z\rangle x$;\\
\item[iv.] $\langle x\times y,x\rangle=0$ and $\langle x\times y,y\rangle=0$;\\
\item[v.] $\langle x\times y,x\times y \rangle=-\langle x,x\rangle \langle y,y\rangle+(\langle x,y\rangle)^{2}$ for all $x, y, z$ in $\mathbb{R}^{3}_{1}$.
\end{enumerate}
We now define "spheres" in $\mathbb{R}^{3}_{1}$ as follows:
{\setlength\arraycolsep{1pt}
\begin{eqnarray}
\mathbb{S}^{2}_{1} & = & \big\{(x_1,x_2,x_3)\in\mathbb{R}^{3}_{1}: x_1^{2}+x_2^{2}-x_3^{2}=1\big\}; \nonumber \\
\mathbb{H}^{2} & = & \big\{(x_1,x_2,x_3)\in\mathbb{R}^{3}_{1}: x_1^{2}+x_2^{2}-x_3^{2}=-1\big\}. \nonumber
\end{eqnarray}}
We call $\mathbb{S}^{2}_{1}$ as de Sitter 2-space and $\mathbb{H}^{2}$ as hyperbolic space.

For a unit speed curve $\alpha$ in $\mathbb{R}^{3}_{1}$, we can define a Frenet frame $\{T(s), N(s), B(s)\}$ associated for each point $s$. Here $T, N$ and $B$ are the tangent, normal and binormal vector fields, respectively. Depending on the causal character of the curve $\alpha$, we have the following Frenet equations, pseudo-spherical Darboux images and curvatures:

Let $\alpha$ be a unit speed space-like curve in $\mathbb{R}^{3}_{1}$. We assume that $T'(s)$ is space-like. The Frenet frame $\{T(s), N(s), B(s)\}$ of $\alpha$ is given by
\begin{eqnarray}
T(s)=\alpha'(s),\ N(s)=\alpha''(s)/\|\alpha''(s)\| \ \textrm{and}\ B(s)=N(s)\times T(s). \nonumber
\end{eqnarray}
We define the curvature of $\alpha$ at $s$ as $\kappa (s)=\|T'(s)\|$ or $\kappa (s)=\langle T'(s), N(s)\rangle$. Also we define the torsion of $\alpha$ at $s$ as $\tau(s)=-\langle N'(s), B(s)\rangle$. Thus the Frenet equations are
\begin{equation}
\left[%
\begin{array}{c}
  T'(s)\\
  N'(s)\\
  B'(s)\\
\end{array}%
\right]
=
\left[%
\begin{array}{ccc}
  0 & \kappa(s) & 0\\
  -\kappa(s) & 0 & \tau(s)\\
  0 & \tau(s) & 0\\
\end{array}%
\right]
\left[%
\begin{array}{c}
  T(s)\\
  N(s)\\
  B(s)\\
\end{array}%
\right]. \nonumber
\end{equation}
The Darboux vector of this space-like curve $\alpha$ is given by $D(s)=-\tau(s)T(s)+\kappa(s)B(s)$. Suppose that $D(s)$ is non-light-like vector for any $s\in I\subset\mathbb{R}$. Then de Sitter Darboux image or pseudo-spherical Darboux image of $\alpha$ \cite{wang} is defined by
\begin{equation}
C: I\subset\mathbb{R} \rightarrow \mathbb{S}^{2}_{1},\ \ s \rightarrow C(s)=\frac{D(s)}{\|D(s)\|}. \nonumber
\end{equation}
For a general parameter $t$ of a space-like space curve $\alpha$, we can find the curvature and torsion as follows:
\begin{equation}
\kappa(t)=\frac{\|\alpha'(t)\times\alpha''(t)\|}{(\langle\alpha'(t),\alpha'(t)\rangle)^{3/2}} \ \textrm{and}\ \tau(t)=\frac{\det(\alpha'(t),\alpha''(t),\alpha'''(t))}{\|\alpha'(t)\times\alpha''(t)\|^2}.
\end{equation}
Also, let $\alpha$ be a unit speed time-like curve in $\mathbb{R}^{3}_{1}$. The Frenet frame $\{T(s), N(s), B(s)\}$ of $\alpha$ is given by
\begin{eqnarray}
T(s)=\alpha'(s),\ N(s)=\alpha''(s)/\|\alpha''(s)\| \ \textrm{and}\ B(s)=T(s)\times N(s). \nonumber
\end{eqnarray}
We define the curvature of $\alpha$ at $s$ as $\kappa (s)=\|T'(s)\|$. Also we define the torsion of $\alpha$ at $s$ as $\tau(s)=\langle N'(s), B(s)\rangle$. Thus the Frenet equations are
\begin{equation}
\left[%
\begin{array}{c}
  T'(s)\\
  N'(s)\\
  B'(s)\\
\end{array}%
\right]
=
\left[%
\begin{array}{ccc}
  0 & \kappa(s) & 0\\
  \kappa(s) & 0 & \tau(s)\\
  0 & -\tau(s) & 0\\
\end{array}%
\right]
\left[%
\begin{array}{c}
  T(s)\\
  N(s)\\
  B(s)\\
\end{array}%
\right]. \nonumber
\end{equation}
The Darboux vector of the time-like curve $\alpha$ is given by $D(s)=\tau(s)T(s)+\kappa(s)B(s)$. We suppose that $D(s)$ is non-light-like vector for any  $s\in I\subset\mathbb{R}$. Then hyperbolic Darboux image or pseudo-spherical Darboux image of $\alpha$ \cite{wang} is defined by
\begin{equation}
C: I\subset\mathbb{R} \rightarrow \mathbb{H}^{2},\ \ s\rightarrow C(s)=\frac{D(s)}{\|D(s)\|}. \nonumber
\end{equation}
For a general parameter $t$ of a time-like space curve $\alpha$, we can find the curvature and the torsion as follows:
\begin{equation}
\kappa(t)=\frac{\|\alpha'(t)\times\alpha''(t)\|}{(-\langle\alpha'(t),\alpha'(t)\rangle)^{3/2}}\ \textrm{and}\ \tau(t)=\frac{\det(\alpha'(t),\alpha''(t),\alpha'''(t))}{\|\alpha'(t)\times\alpha''(t)\|^2}.
\end{equation}
In  $\mathbb{R}^{3}_{1}$, we say that if a space-like or time-like curve $\alpha$ is a helix, then $\tau/\kappa$ is a constant function. Conversely, let $\alpha$ be a space-like or time-like curve with non-light-like normal vector. If $\tau/\kappa$ is constant, then $\alpha$ is a helix. In addition, a space-like or time-like curve $\alpha$ is a Bertrand curve if and only if there are non-zero real constants $A, B$ such that $A\kappa(s)+B\tau(s)=1$ for any $s\in I\subset\mathbb{R}$.


\section{Space-like constant slope surfaces lying in the space-like cone and space-like Bertrand curves in Minkowski 3-space}

In this section, we define the notions of Lorentzian Sabban frames and de Sitter evolutes of unit speed space-like curves on de Sitter 2-space $\mathbb{S}^{2}_{1}$. Also we investigate the invariants and geometric properties of these curves. We show that space-like Bertrand curves can be constructed from unit speed space-like curves on $\mathbb{S}^{2}_{1}$. Moreover, we give a relation between space-like Bertrand curves and helices. We show that de Sitter Darboux images of Bertrand curves are equal to de Sitter evolutes of unit speed space-like curves on $\mathbb{S}^{2}_{1}$. Finally, we obtain some relations between space-like Bertrand curves and space-like constant slope surfaces lying in the space-like cone.

Now we define a pseudo-orthonormal frame along a space-like curve on $\mathbb{S}^{2}_{1}$. Let  $f:I\rightarrow\mathbb{S}^{2}_{1}$ be unit speed space-like curve with the tangent vector $t(v)=f'(v)$. Here, $v$ is arc-length parameter of $f$. We now set a vector $s(v)=f(v)\times t(v)$ and as a consequence   $s(v)\times t(v)=f(v)$, where $f$ denotes the position vector of the curve. By definition of the space-like curve $f$, we have a Lorentzian Sabban frame    $\{f(v), t(v), s(v)\}$ along $f$. Thus we have the following pseudo-spherical Frenet-Serret formulae of $f:$
\begin{equation}
\left[%
\begin{array}{c}
  f'(v)\\
  t'(v)\\
  s'(v)\\
\end{array}%
\right]
=
\left[%
\begin{array}{ccc}
  0 & 1 & 0\\
  -1 & 0 & -\kappa_{g}(v)\\
  0 & -\kappa_{g}(v) & 0\\
\end{array}%
\right]
\left[%
\begin{array}{c}
  f(v)\\
  t(v)\\
  s(v)\\
\end{array}%
\right],
\end{equation}
where $\kappa_{g}(v)$ is the geodesic curvature of the curve $f$ on $\mathbb{S}^{2}_{1}$ which is given by $\kappa_{g}(v)=\det(f(v), t(v), t'(v))$.

Under the assumption that $\kappa_{g}^{2}(v)>1$, we define a curve on $\mathbb{S}^{2}_{1}$ as
\begin{equation}
d_{f}(v)=\frac{-\kappa_{g}(v)f(v)-s(v)}{\sqrt{\kappa_{g}^{2}(v)-1}}. \nonumber
\end{equation}
We call $d_{f}$ as de Sitter evolute or pseudo-spherical evolute of the space-like curve $f$.

{\bf 3.1. Space-like height function of unit speed space-like curves on $\mathbb{S}^{2}_{1}$}

In this section, we define a function on a space-like curve $f:I\rightarrow\mathbb{S}^{2}_{1}$ by using the similar methods in \cite{izumiya2}. Now we define the function $H^{s}:I\times\mathbb{S}^{2}_{1}\rightarrow\mathbb{R}$ by $H^{s}(v,u)=\langle f(v),u\rangle$. We call $H^{s}$ as the space-like height function of $f$ and we denote $(h_{u}^{s})(v)=H^{s}(v,u)$.

We now give the following proposition.
\begin{proposition}
Let $f:I\rightarrow\mathbb{S}^{2}_{1}$ be a unit speed space-like curve. For any $(v,u)\in I\times\mathbb{S}^{2}_{1}$:
\begin{enumerate}
\item[a.] $(h_{u}^{s})'(v)=0 \ \textrm{if and only if}\ u\in span\{f(v), s(v)\}$;\\
\item[b.] $(h_{u}^{s})'(v)=(h_{u}^{s})''(v)=0 \ \textrm{if and only if}\ u=\pm\frac{\kappa_{g}(v)f(v)+s(v)}{\sqrt{\kappa_{g}^{2}(v)-1}}, \ \textrm{where}\ \kappa_{g}^{2}(v)>1$. \nonumber
\end{enumerate}
\end{proposition}
\begin{proof}
Using the pseudo-spherical Frenet formulae in (5), we obtain
\begin{enumerate}
\item[i.] $(h_{u}^{s})'(v)=\langle t(v),u\rangle$;\\
\item[ii.] $(h_{u}^{s})''(v)=\langle -f(v)-\kappa_{g}(v)s(v),u\rangle.$ \nonumber
\end{enumerate}
The case (a) can be obtained from Eq. (i). By using the case (a), there exists $\lambda, \mu\in \mathbb{R}$ such that $u=\lambda f(v)+\mu s(v)$. From Eq. (ii), we get
{\setlength\arraycolsep{1pt}
\begin{eqnarray}
0 & = & \langle-f(v)-\kappa_{g}(v)s(v), \lambda f(v)+\mu s(v)\rangle \nonumber \\
& = & -\lambda \langle f(v), f(v)\rangle-\mu \kappa_{g}(v)\langle s(v),s(v)\rangle\nonumber \\
& = & -\lambda+\mu\kappa_{g}(v). \nonumber
\end{eqnarray}}
Hence we have $u=\mu(\kappa_{g}(v)f(v)+s(v))$. Since $\langle u,u\rangle=1$, we get
\begin{equation}
\mu=\pm\frac{1}{\sqrt{\kappa_{g}^{2}(v)-1}}. \nonumber
\end{equation}
Thus, we obtain
\begin{equation}
u=\pm\frac{1}{\sqrt{\kappa_{g}^{2}(v)-1}}(\kappa_{g}(v)f(v)+s(v)). \nonumber
\end{equation}
The proof is completed.
\end{proof}

{\bf 3.2. Spherical invariants of unit speed space-like curves on de Sitter 2-space $\mathbb{S}^{2}_{1}$}

In this section, we investigate the geometric properties of de Sitter evolutes of the unit speed space-like curves on $\mathbb{S}^{2}_{1}$ by using the similar methods in \cite{izumiya2}.

For any $r\in\mathbb{R}$ and $u_{0}\in\mathbb{S}^{2}_{1}$, we denote $PS^{1}(u_{0},r)=\{u\in\mathbb{S}^{2}_{1}: \langle u,u_{0}\rangle=r\}$. We call $PS^{1}(u_{0},r)$ as a pseudo-circle whose center is $u_{0}$ on $\mathbb{S}^{2}_{1}$.

Thus we state the following proposition.

\begin{proposition}
Let $f:I\rightarrow\mathbb{S}^{2}_{1}$ be a unit speed space-like curve with $\kappa_{g}^{2}(v)>1$. Then $\kappa_{g}'(v)\equiv0$ if and only if $u_{0}=\pm(\kappa_{g}(v)f(v)+s(v))/\sqrt{\kappa_{g}^{2}(v)-1}$ are constant vectors. Under this condition, $f$ is a part of a pseudo-circle whose center is $u_{0}$ on $\mathbb{S}^{2}_{1}$.
\end{proposition}
\begin{proof}
Let us denote
\begin{equation}
P_{\pm}(v)=\pm u_{0}=\pm\frac{1}{\sqrt{\kappa_{g}^{2}(v)-1}}(\kappa_{g}(v)f(v)+s(v)). \nonumber
\end{equation}
If we take the derivative of this equation with respect to $v$, then we have
\begin{equation}
P_{\pm}'(v)=\pm\kappa_{g}'(v)\frac{(f(v)+\kappa_{g}(v)s(v))}{(\kappa_{g}^{2}(v)-1)^{3/2}}. \nonumber
\end{equation}
Therefore, $P_{\pm}'(v)\equiv0$ if and only if $\kappa_{g}'(v)\equiv0$.

Under this condition, if we take
\begin{equation}
r=\pm\frac{\kappa_{g}(v)}{\sqrt{\kappa_{g}^{2}(v)-1}} \ \textrm{and}\ u_{0}=\pm\frac{\kappa_{g}(v)f(v)+s(v)}{\sqrt{\kappa_{g}^{2}(v)-1}}, \nonumber
\end{equation}
then $f(v)$ is a part of the pseudo-circle $PS^{1}(u_{0},r)$. This completes the proof.
\end{proof}

Let $f:I\rightarrow\mathbb{S}^{2}_{1}$ be a unit space-like curve with $\kappa_{g}^{2}(v)>1$. For any $v_{0}\in I$, let us consider the pseudo-circle $PS^{1}(u_{0},r_{0}^{\pm})$, where $u_{0}=d_{f}(v_{0})$ and $r_{0}=\kappa_{g}(v_{0})/\sqrt{\kappa_{g}^{2}(v_{0})-1}$.

Then we have the following proposition.

\begin{proposition}
Under the above notations, $f$ and $PS^{1}(u_{0},r_{0})$ have at least a 3-point contact at $f(v_{0})$.
\end{proposition}
\begin{proof}
In Proposition 1, the case (b) says that $f$ and $PS^{1}(u_{0},r_{0})$ have at least a 3-point contact at $f(v_{0})$. This completes the proof.
\end{proof}
\begin{remark}
We call $PS^{1}(u_{0},r_{0})$ in Proposition 3 as the pseudo-circle of geodesic curvature and its center $u_{0}$ is called as the center of geodesic curvature. As a result, de Sitter evolute $d_{f}(v_{0})$ is the locus of the center of geodesic curvature.
\end{remark}

{\bf 3.3. Space-like constant slope surfaces lying in the space-like cone and space-like Bertrand curves}

Firstly, we can express the following lemma.

\begin{lemma}
\rm Let $f:I\rightarrow \mathbb{S}^{2}_{1}$ be a unit speed space-like curve. Then
\begin{eqnarray} \label{eq5}
\tilde{\gamma}(v)=a\int_{0}^{v}f(t)dt+a\tanh\xi_{1}\int_{0}^{v}f(t)\times f'(t)dt
\end{eqnarray}
is a space-like Bertrand curve, where $a,\ \xi_{1}=\xi_{1}(u)=\tanh\theta \ln u$ are constants and $\theta$ is a constant different from $0$. Moreover, all space-like Bertrand curves can be constructed by this method.
\end{lemma}

\begin{proof}
We now find the curvature and torsion of $\tilde{\gamma}(v)$. If we take the derivative of (6) three times with respect to $v$, we have
{\setlength\arraycolsep{1pt}
\begin{eqnarray}
\tilde{\gamma}'(v) & = & a(f(v)+\tanh\xi_{1}s(v)); \nonumber \\
\tilde{\gamma}''(v) & = & a(1-\tanh\xi_{1}\kappa_{g}(v))t(v); \nonumber \\
\tilde{\gamma}'''(v) & = & -a(1-\tanh\xi_{1}\kappa_{g}(v))f(v)-a\tanh\xi_{1}\kappa_{g}'(v)t(v)-a(\kappa_{g}(v)-\tanh\xi_{1}\kappa_{g}^{2}(v))s(v). \nonumber
\end{eqnarray}}
Therefore, by using (3), $\kappa(v)$ and $\tau(v)$ are found as the following:
\begin{equation}
\kappa(v)=\varepsilon\frac{\cosh^{2}\xi_{1}(1-\tanh\xi_{1}\kappa_{g}(v))}{a} \ \textrm{and}\ \tau(v)=\frac{\cosh^{2}\xi_{1}(\kappa_{g}(v)-\tanh\xi_{1})}{a},
\end{equation}
where $\varepsilon=\pm1$. It follows that $a(\varepsilon\kappa(v)+\tanh\xi_{1}\tau(v))=1$, so $\tilde{\gamma}(v)$ is a space-like Bertrand curve.

Conversely, let $\tilde{\gamma}(s)$ be a space-like Bertrand curve. Then, there are real constants $A, B$ different from $0$ such that $A\kappa(s)+B\tau(s)=1$. In this equation, we put $A=a$ and $B=a\tanh\xi_{1}$. Assume that $a>0$ and choose $\varepsilon=\pm1$ with $\varepsilon\cosh\xi_{1}/a>0$.

We now consider the Frenet frame $\{T(s), N(s), B(s)\}$ for the space-like curve $\tilde{\gamma}(s)$. In this trihedron $T(s), N(s)$ are space-like vectors and $B(s)$ is a time-like vector. Thus we can write
\begin{equation}
T(s)\times N(s)=-B(s) \ \textrm{and}\ B(s)\times N(s)=-T(s). \nonumber
\end{equation}
Now we define a space-like curve on $\mathbb{S}^{2}_{1}$ as
\begin{equation}
f(s)=\varepsilon(\cosh\xi_{1}T(s)+\sinh\xi_{1}B(s)). \nonumber
\end{equation}
So we have
\begin{equation}
f'(s)=\varepsilon\cosh\xi_{1}(\kappa(s)+\tanh\xi_{1}\tau(s))N(s)=\frac{\varepsilon}{a}\cosh\xi_{1}N(s). \nonumber
\end{equation}
Let $v$ be the arc-length parameter of $f$. Then we have $dv/ds=\varepsilon\cosh\xi_{1}/a$. Also we have
\begin{equation}
af(s)\frac{dv}{ds}=\cosh\xi_{1}(\cosh\xi_{1}T(s)+\sinh\xi_{1}B(s))
\end{equation}
and
{\setlength\arraycolsep{1pt}
\begin{eqnarray}
a \tanh\xi_{1}f(s)\times\frac{df}{dv} \frac{dv}{ds} & = & a\tanh\xi_{1}\varepsilon(\cosh\xi_{1}T(s)+\sinh\xi_{1}B(s))\times\frac{\varepsilon}{a}\cosh\xi_{1}N(s) \nonumber \\
& = & \sinh\xi_{1}(-\cosh\xi_{1}B(s)-\sinh\xi_{1}T(s)).
\end{eqnarray}}
By using (8) and (9), we conclude that
{\setlength\arraycolsep{1pt}
\begin{eqnarray}
a\int_{0}^{v}f(t)dt+a\tanh\xi_{1}\int_{0}^{v}f(t)\times f'(t)dt & = &  \int_{s_{0}}^{s}\cosh\xi_{1}(\cosh\xi_{1}T(t)+\sinh\xi_{1}B(t)) dt \nonumber \\
&& {}+\int_{s_{0}}^{s}\sinh\xi_{1}(-\cosh\xi_{1}B(t)-\sinh\xi_{1}T(t))dt \nonumber \\
& = & \int_{s_{0}}^{s}T(t)dt=\tilde{\gamma}(s). \nonumber
\end{eqnarray}}
This completes the proof.
\end{proof}
As a consequence of this lemma, we can give a relation between space-like Bertrand curves and helices.
\begin{corollary}
The unit speed space-like curve $f$ on $\mathbb{S}^{2}_{1}$ is a part of a pseudo-circle if and only if the corresponding space-like Bertrand curve is a helix.
\end{corollary}
\begin{proof}
By using (7), we have
\begin{equation}
\kappa'(v)=-\varepsilon\frac{\sinh2\xi_{1}\kappa_{g}'(v)}{2a} \ \textrm{and}\ \tau'(v)=\frac{\cosh^{2}\xi_{1}\kappa_{g}'(v)}{a}. \nonumber
\end{equation}
From Proposition 2, the unit speed space-like curve $f$ on $\mathbb{S}^{2}_{1}$ is a part of a pseudo-circle if and only if $\kappa'(v)\equiv0$. This condition is equivalent to the condition that both $\kappa(v)$ and $\tau(v)$ are non-zero constants. This completes the proof.
\end{proof}
Now we give the following proposition.
\begin{proposition}
Let $f:I\rightarrow\mathbb{S}^{2}_{1}$ be a unit speed space-like curve and $\tilde{\gamma}:I\rightarrow\mathbb{R}^{3}_{1}$ be a space-like Bertrand curve corresponding to $f$. Then de Sitter Darboux image of $\tilde{\gamma}$ is equal to de Sitter evolute of $f$.
\end{proposition}
\begin{proof}
By (7), we have
\begin{equation}
\kappa(v)=\varepsilon\frac{\cosh^{2}\xi_{1}(1-\tanh\xi_{1}\kappa_{g}(v))}{a} \ \textrm{and}\ \tau(v)=\frac{\cosh^{2}\xi_{1}(\kappa_{g}(v)-\tanh\xi_{1})}{a}. \nonumber
\end{equation}
For the space-like curve $\tilde{\gamma}$, we have
\begin{equation}
T(v)=a(f(v)+\tanh\xi_{1}s(v))\frac{dv}{ds} \ \textrm{and}\  N(v)=\varepsilon t(v). \nonumber
\end{equation}
Also we get
\begin{equation}
B(v)=N(v)\times T(v)=-\varepsilon a\frac{dv}{ds}(s(v)+\tanh\xi_{1}f(v)). \nonumber
\end{equation}
Thus we can show that
\begin{equation}
D(v)=-\tau(v)T(v)+\kappa(v)B(v)=\frac{dv}{ds}(-\kappa_{g}(v)f(v)-s(v)). \nonumber
\end{equation}	
As a result, we have $C(v)=D(v)/\|D(v)\|=d_{f}(v)$. This completes the proof.
\end{proof}	
We now state the relations between space-like Bertrand curves and space-like constant slope space-like surfaces lying in the space-like cone.
\begin{theorem}
Let $f:I\rightarrow\mathbb{S}^{2}_{1}$ be a unit speed space-like curve and $\tilde{\gamma}:I\rightarrow\mathbb{R}^{3}_{1}$ be a space-like Bertrand curve corresponding to $f$. Then $\tilde{\gamma}'(v)$ lies on the space-like constant slope surface $x(u,v)$ lying in the space-like cone.
\end{theorem}
\begin{proof}
By Lemma 1, if we take the derivative of (6) with respect to $v$, we obtain
\begin{equation}
\tilde{\gamma}'(v)=af(v)+a\tanh\xi_{1}f(v)\times f'(v). \nonumber
\end{equation}
In this equation, we can take $a=u\cosh\theta\cosh\xi_{1}$ and thus $a\tanh\xi_{1}=u\cosh\theta\sinh\xi_{1}$, where $u, \theta$ are constants. Thus, from (1), $\tilde{\gamma}'(v)$ is $v$-parameter curve of space-like constant slope surface $x(u,v)$ lying in the space-like cone and $\tilde{\gamma}'(v)$ lies on it. The proof is completed.
\end{proof}
Also the following theorem can be given.
\begin{theorem}
Let $x:S\rightarrow\mathbb{R}^{3}_{1}$ be a space-like constant slope surface immersed in $\mathbb{R}^{3}_{1}$ which lies in the space-like cone. If $x(v)$ is $v$-parameter curve of space-like constant slope surface $x(u,v)$ lying in the space-like cone, then $\int_{0}^{v}x(v)dv$ is a space-like Bertrand curve.
\end{theorem}
\begin{proof}
From (1), we get
\begin{equation}
x(v)=u\cosh\theta\cosh\xi_{1}f(v)+u\cosh\theta\sinh\xi_{1}f(v)\times f'(v). \nonumber
\end{equation}
for $u$=constant, where $\xi_{1}=\xi_{1}(u)=\tanh \theta \ln u$. By integrating $x(v)$, we have
\begin{equation}
\int_{0}^{v}x(v)dv=u\cosh\theta\cosh\xi_{1}\int_{0}^{v}f(v)dv+u\cosh\theta\sinh\xi_{1}\int_{0}^{v}f(v)\times f'(v)dv. \nonumber
\end{equation}
Since the coefficients of $f(v)$ and $f(v)\times f'(v)$ are constants, here we can take $u\cosh\theta\cosh\xi_{1}=a$ and so $u\cosh\theta\sinh\xi_{1}=a\tanh\xi_{1}$. Therefore we obtain
\begin{equation}
\int_{0}^{v}x(v)dv=a\int_{0}^{v}f(v)dv+a\tanh\xi_{1}\int_{0}^{v}f(v)\times f'(v)dv. \nonumber
\end{equation}
By Lemma 1, $\int_{0}^{v}x(v)dv$ is a space-like Bertrand curve. This completes the proof.
\end{proof}
We now give an example for space-like constant slope surfaces and space-like Bertrand curves and draw their pictures by using Mathematica.
\begin{example}
Let us take the unit speed space-like curve $f(v)=(\sin v, \cos v, 0)$ on $\mathbb{S}^{2}_{1}$. Then we have $f(v)\times f'(v)=(0,0,1)$. By using (1), the space-like constant slope surface lying in the space-like cone is given by
\begin{equation}
x(u,v)=u\cosh\theta \big(\cosh(\tanh\theta\ln u)\sin v,\cosh(\tanh\theta\ln u)\cos v, \sinh(\tanh\theta\ln u) \big). \nonumber
\end{equation}
For $\theta=1.5$, the picture of this surface is drawn by Figure~\ref{fig1}.
\begin{figure}[htbp]
\includegraphics[height=60mm]{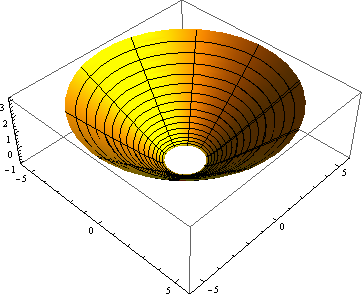}
\caption {The space-like constant slope surface lying in the space-like cone}
\label{fig1}
\end{figure}
Also, for $u=e$, the space-like Bertrand curve is given by
\begin{equation}
\int_{0}^{v}x(v)dv=e\cosh(1.5) \big(-\cosh(\tanh(1.5))(\cos v-1),\cosh(\tanh(1.5))\sin v, \sinh(\tanh(1.5))v \big). \nonumber
\end{equation}
Since the space-like curve $f(v)$ is a part of a pseudo-circle on $\mathbb{S}^{2}_{1}$, from Corollary 1, this space-like Bertrand curve is a helix. The picture of this curve is drawn by Figure~\ref{fig2}.
\begin{figure}[htbp]
  \includegraphics[height=60mm]{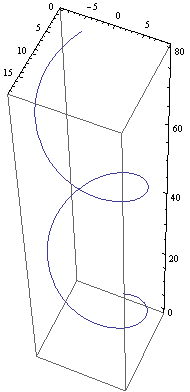}
\caption{The space-like Bertrand curve}
\label{fig2}
\end{figure}

\end{example}


\section{Space-like constant slope surfaces lying in the time-like cone and time-like Bertrand curves in Minkowski 3-space}

As in Section 3, firstly, we recall the notions of Lorentzian Sabban frames and hyperbolic evolutes of unit speed space-like curves on the hyperbolic space $\mathbb{H}^{2}$. Furthermore, we show that time-like Bertrand curves can be constructed from unit speed space-like curves on $\mathbb{H}^{2}$. Also we give a relation between time-like Bertrand curves and helices and show that hyperbolic Darboux images of Bertrand curves are equal to hyperbolic evolutes of the unit speed space-like curves on $\mathbb{H}^{2}$. Finally, we obtain the relations between time-like Bertrand curves and space-like constant slope surfaces lying in the time-like cone.

Let $g:I \rightarrow \mathbb{H}^{2}$ be a unit speed space-like curve. We denote $v$ as the arc-length parameter of $g$. Let us denote $t(v)=g'(v)$ and we call $t(v)$ the unit tangent vector of $g$ at $v$. We now set a vector $s(v)=g(v)\times t(v)$ and as a consequence $s(v)\times t(v)=g(v)$, where $g$ denotes the position vector of the curve. By definition of the space-like curve $g$, we have a Lorentzian Sabban frame $\{g(v), t(v), s(v)\}$ along $g$. Then we have the following pseudo-spherical Frenet-Serret formulae of $g$:
\begin{equation}
\left[%
\begin{array}{c}
  g'(v)\\
  t'(v)\\
  s'(v)\\
\end{array}%
\right]
=
\left[%
\begin{array}{ccc}
  0 & 1 & 0\\
  1 & 0 & \kappa_{g}(v)\\
  0 & -\kappa_{g}(v) & 0\\
\end{array}%
\right]
\left[%
\begin{array}{c}
  g(v)\\
  t(v)\\
  s(v)\\
\end{array}%
\right],
\end{equation}
where $\kappa_{g}(v)$ is the geodesic curvature of the curve $g$ on $\mathbb{H}^{2}$ which is given by $\kappa_{g}(v)=\det(g(v), t(v), t'(v))$ \cite{izumiya2}.

Under the assumption that $\kappa_{g}^{2}(v)>1$, we define a curve on $\mathbb{H}^{2}$ as
\begin{equation}
h_{g}(v)=\frac{\kappa_{g}(v)g(v)+s(v)}{\sqrt{\kappa_{g}^{2}(v)-1}}. \nonumber
\end{equation}
We call $h_{g}$ as hyperbolic evolute or pseudo-spherical evolute of the space-like curve $g$ \cite{izumiya2}.

For any $r\in\mathbb{R}$ and $u_{0}\in\mathbb{H}^{2}$, we denote $PS^{1}(u_{0},r)=\{u\in\mathbb{H}^{2}: \langle u,u_{0}\rangle=r\}$. We call $PS^{1}(u_{0},r)$ as a pseudo-circle whose center is $u_{0}$ on $\mathbb{H}^{2}$ \cite{izumiya2}.

Let $g:I\rightarrow\mathbb{H}^{2}$ be a unit space-like curve with $\kappa_{g}^{2}(v)>1$. For any $v_{0}\in I$, we consider the pseudo-circle $PS^{1}(u_{0},r_{0}^{\pm})$, where $u_{0}=h_{g}(v_{0})$ and $r_{0}=-\kappa_{g}(v_{0})/\sqrt{\kappa_{g}^{2}(v_{0})-1}$. We call $PS^{1}(u_{0},r_{0})$ as the pseudo-circle of geodesic curvature and its center $u_{0}$ is called as the center of geodesic curvature. Thus hyperbolic evolute $h_{g}(v_{0})$ is the locus of the center of geodesic curvature (see \cite{izumiya2}).

Now we give the following lemma.
\begin{lemma}
\rm Let $g:I\rightarrow \mathbb{H}^{2}$ be a unit speed space-like curve. Then
\begin{eqnarray} \label{eq5}
\tilde{\gamma}(v)=a\int_{0}^{v}g(t)dt+a\tanh\xi_{2}\int_{0}^{v}g(t)\times g'(t)dt
\end{eqnarray}
is a time-like Bertrand curve, where $a,\ \xi_{2}=\xi_{2}(u)=\coth\theta \ln u$ are constants and $\theta$ is a constant different from $0$. Moreover, all time-like Bertrand curves can be constructed by this method.
\end{lemma}
\begin{proof}
We now find the curvature and torsion of $\tilde{\gamma}(v)$. If we take the derivative of (11) three times with respect to $v$, then we have
{\setlength\arraycolsep{1pt}
\begin{eqnarray}
\tilde{\gamma}'(v) & = & a(g(v)+\tanh\xi_{2}s(v)); \nonumber \\
\tilde{\gamma}''(v) & = & a(1-\tanh\xi_{2}\kappa_{g}(v))t(v); \nonumber \\
\tilde{\gamma}'''(v) & = & a(1-\tanh\xi_{2}\kappa_{g}(v))g(v)-a\tanh\xi_{2}\kappa_{g}'(v)t(v)+a(\kappa_{g}(v)-\tanh\xi_{2}\kappa_{g}^{2}(v))s(v). \nonumber
\end{eqnarray}}
Therefore, by using (4), we can find $\kappa(v)$ and $\tau(v)$ as follows:
\begin{equation}
\kappa(v)=\varepsilon\frac{\cosh^{2}\xi_{2}(1-\tanh\xi_{2}\kappa_{g}(v))}{a} \ \textrm{and}\ \tau(v)=\frac{\cosh^{2}\xi_{2}(\kappa_{g}(v)-\tanh\xi_{2})}{a},
\end{equation}
where $\varepsilon=\pm1$. It follows that $a(\varepsilon\kappa(v)+\tanh\xi_{2}\tau(v))=1$, so $\tilde{\gamma}(v)$ is a time-like Bertrand curve.

Conversely, let $\tilde{\gamma}(s)$ be a time-like Bertrand curve. Then, there are real constants $A, B$ different from $0$ such that $A\kappa(s)+B\tau(s)=1$. In this equation, we put $A=a$ and $B=a\tanh\xi_{2}$. Assume that $a>0$ and choose $\varepsilon=\pm1$ with $\varepsilon\cosh\xi_{2}/a>0$.

Let us consider the Frenet frame $\{T(s), N(s), B(s)\}$ for the time-like curve $\tilde{\gamma}(s)$. In this trihedron $T(s)$ is a time-like vector, $N(s)$ and $B(s)$ are space-like vectors. Thus we can write
\begin{equation}
T(s)\times N(s)=B(s) \ \textrm{and}\ B(s)\times N(s)=T(s). \nonumber
\end{equation}
Now we define a space-like curve on $\mathbb{H}^{2}$ as
\begin{equation}
g(s)=\varepsilon(\cosh\xi_{2}T(s)-\sinh\xi_{2}B(s)). \nonumber
\end{equation}
Thus we have
\begin{equation}
g'(s)=\varepsilon\cosh\xi_{2}(\kappa(s)+\tanh\xi_{2}\tau(s))N(s)=\frac{\varepsilon}{a}\cosh\xi_{2}N(s). \nonumber
\end{equation}
Let $v$ be the arc-length parameter of $g$, thus we have $dv/ds=\varepsilon\cosh\xi_{2}/a$. Moreover we have
\begin{equation}
ag(s)\frac{dv}{ds}=\cosh\xi_{2}(\cosh\xi_{2}T(s)-\sinh\xi_{2}B(s))
\end{equation}
and
{\setlength\arraycolsep{1pt}
\begin{eqnarray}
a \tanh\xi_{2}g(s)\times\ \frac{dg}{dv}\frac{dv}{ds} & = & a\tanh\xi_{2}\varepsilon(\cosh\xi_{2}T(s)-\sinh\xi_{2}B(s))\times\frac{\varepsilon}{a}\cosh\xi_{2}N(s) \nonumber \\
& = & \sinh\xi_{2}(\cosh\xi_{2}B(s)-\sinh\xi_{2}T(s)).
\end{eqnarray}}
By using (13) and (14), we conclude that
{\setlength\arraycolsep{1pt}
\begin{eqnarray}
a\int_{0}^{v}g(t)dt+a\tanh\xi_{2}\int_{0}^{v}g(t)\times g'(t)dt & = &  \int_{s_{0}}^{s}\cosh\xi_{2}(\cosh\xi_{2}T(t)-\sinh\xi_{2}B(t)) dt \nonumber \\
&& {}+\int_{s_{0}}^{s}\sinh\xi_{2}(\cosh\xi_{2}B(t)-\sinh\xi_{2}T(t))dt \nonumber \\
& = & \int_{s_{0}}^{s}T(t)dt=\tilde{\gamma}(s). \nonumber
\end{eqnarray}}
This completes the proof.
\end{proof}
As a consequence of this lemma, we can give a relation between time-like Bertrand curves and helices as follows.
\begin{corollary}
The unit speed space-like curve $g$ on $\mathbb{H}^{2}$ is a part of a pseudo-circle if and only if the corresponding time-like Bertrand curve is a helix.
\end{corollary}
\begin{proof}
By using (12), we get
\begin{equation}
\kappa'(v)=-\varepsilon\frac{\sinh2\xi_{2}\kappa_{g}'(v)}{2a} \ \textrm{and}\ \tau'(v)=\frac{\cosh^{2}\xi_{2}\kappa_{g}'(v)}{a}. \nonumber
\end{equation}
The unit speed space-like curve $g$ on $\mathbb{H}^{2}$ is a part of a pseudo-circle if and only if $\kappa'(v)\equiv0$ \cite{izumiya2}. This condition is equivalent to the condition that both $\kappa(v)$ and $\tau(v)$ are non-zero constants. This completes the proof.
\end{proof}
Also we have the following proposition.
\begin{proposition}
Let $g:I\rightarrow\mathbb{H}^{2}$ be a unit speed space-like curve and $\tilde{\gamma}:I\rightarrow\mathbb{R}^{3}_{1}$ be a time-like Bertrand curve corresponding to $g$. Then hyperbolic Darboux image of $\tilde{\gamma}$ is equal to hyperbolic evolute of $g$.
\end{proposition}
\begin{proof}
From (12), we have
\begin{equation}
\kappa(v)=\varepsilon\frac{\cosh^{2}\xi_{2}(1-\tanh\xi_{2}\kappa_{g}(v))}{a} \ \textrm{and}\ \tau(v)=\frac{\cosh^{2}\xi_{2}(\kappa_{g}(v)-\tanh\xi_{2})}{a}. \nonumber
\end{equation}
For the time-like curve $\tilde{\gamma}$, we obtain
\begin{equation}
T(v)=a(g(v)+\tanh\xi_{2}s(v))\frac{dv}{ds} \ \textrm{and}\  N(v)=\varepsilon t(v). \nonumber
\end{equation}
From here, we get
\begin{equation}
B(v)=T(v)\times N(v)=\varepsilon a\frac{dv}{ds}(s(v)+\tanh\xi_{2}g(v)). \nonumber
\end{equation}
Thus we can easily obtain
\begin{equation}
D(v)=\tau(v)T(v)+\kappa(v)B(v)=\frac{dv}{ds}(\kappa_{g}(v)g(v)+s(v)). \nonumber
\end{equation}	
As a result, we have $C(v)=D(v)/\|D(v)\|=h_{g}(v)$. This completes the proof.
\end{proof}	
We now state the relations between time-like Bertrand curves and space-like constant slope space-like surfaces lying in the time-like cone.
\begin{theorem}
Let $g:I\rightarrow\mathbb{H}^{2}$ be a unit speed space-like curve and $\tilde{\gamma}:I\rightarrow\mathbb{R}^{3}_{1}$ be a time-like Bertrand curve corresponding to $g$. Then $\tilde{\gamma}'(v)$ lies on the space-like constant slope surface $x(u,v)$ lying in the time-like cone.
\end{theorem}
\begin{proof}
By Lemma 2, if we take the derivative of (11) with respect to $v$, then we have
\begin{equation}
\tilde{\gamma}'(v)=ag(v)+a\tanh\xi_{2}g(v)\times g'(v). \nonumber
\end{equation}
In this equation, we can take $a=u\sinh\theta\cosh\xi_{2}$ and so $a\tanh\xi_{2}=u\sinh\theta\sinh\xi_{2}$, where $u, \theta$ are constants. Thus, from (2), $\tilde{\gamma}'(v)$ is $v$-parameter curve of space-like constant slope surface $x(u,v)$ lying in the time-like cone and $\tilde{\gamma}'(v)$ lies on it. This completes the proof.
\end{proof}
Also we have the following theorem.
\begin{theorem}
Let $x:S\rightarrow\mathbb{R}^{3}_{1}$ be a space-like constant slope surface immersed in $\mathbb{R}^{3}_{1}$ which lies in the time-like cone. If $x(v)$ is $v$-parameter curve of space-like constant slope surface $x(u,v)$ lying in the time-like cone, then $\int_{0}^{v}x(v)dv$ is a time-like Bertrand curve.
\end{theorem}
\begin{proof}
From (2), we get
\begin{equation}
x(v)=u\sinh\theta\cosh\xi_{2}g(v)+u\sinh\theta\sinh\xi_{2}g(v)\times g'(v). \nonumber
\end{equation}
for $u$=constant, where $\xi_{2}=\xi_{2}(u)=\coth \theta \ln u$. By integrating $x(v)$, we have
\begin{equation}
\int_{0}^{v}x(v)dv=u\sinh\theta\cosh\xi_{2}\int_{0}^{v}g(v)dv+u\sinh\theta\sinh\xi_{2}\int_{0}^{v}g(v)\times g'(v)dv. \nonumber
\end{equation}
Since the coefficients of $g(v)$ and $g(v)\times g'(v)$ are constants, here we can take $u\sinh\theta\cosh\xi_{2}=a$ and so $u\sinh\theta\sinh\xi_{2}=a\tanh\xi_{2}$. Therefore we obtain
\begin{equation}
\int_{0}^{v}x(v)dv=a\int_{0}^{v}g(v)dv+a\tanh\xi_{2}\int_{0}^{v}g(v)\times g'(v)dv. \nonumber
\end{equation}
By Lemma 2, $\int_{0}^{v}x(v)dv$ is a time-like Bertrand curve. This completes the proof.
\end{proof}
Now we give the following example.
\begin{example}
Let us take the unit speed space-like curve as $g(v)=(\sinh v, 0, \cosh v)$ on $\mathbb{H}^{2}$. Then we have $g(v)\times g'(v)=(0,1,0)$. By using (2), the space-like constant slope surface lying in the time-like cone is given by
\begin{equation}
x(u,v)=u\sinh\theta \big(\cosh(\coth\theta\ln u)\sinh v, \sinh(\coth\theta\ln u), \cosh(\coth\theta\ln u)\cosh v\big). \nonumber
\end{equation}
For $\theta=1.5$, the picture of this surface is drawn by Figure~\ref{fig3}.
\begin{figure}[htbp]
\includegraphics[height=55mm]{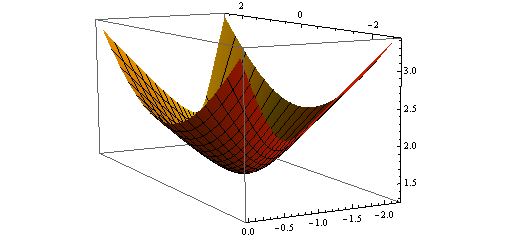}
\caption {The space-like constant slope surface lying in the time-like cone}
\label{fig3}
\end{figure}
Also, for $u=e$, the time-like Bertrand curve is given by
\begin{equation}
\int_{0}^{v}x(v)dv=e\sinh(1.5) \big(\cosh(\coth(1.5))(\cosh v-1), \sinh(\coth(1.5))v, \cosh(\coth(1.5))\sinh v\big). \nonumber
\end{equation}
Since the space-like curve $g(v)$ is a part of a pseudo-circle on $\mathbb{H}^{2}$, from Corollary 2, this time-like Bertrand curve is a helix. The picture of this curve is drawn by Figure~\ref{fig4}.
\begin{figure}[htbp]
  \includegraphics[height=70mm]{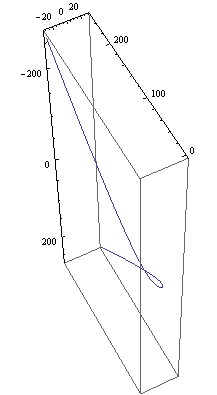}
\caption{The time-like Bertrand curve}
\label{fig4}
\end{figure}

\end{example}

\newpage

{\bf Acknowledgement.} This paper is a part of the doctoral thesis of the first author.

\end{document}